\documentclass[runningheads]{llncs}

\usepackage{graphicx}

\usepackage{wrapfig}

\usepackage{amsmath,amssymb}
\usepackage{tikz}
\usepackage{todonotes}

\newcommand{\freq}{\omega}
\newcommand{\acc}{\nu}
\newcommand{\slope}{\acc}

\DeclareMathOperator{\diff}{d}

\newcommand{\R}{\mathbb{R}}
\newcommand{\C}{\mathbb{C}}

\newcommand{\abs}[1]{\left\lvert#1\right\rvert}
\newcommand{\norm}[1]{\left\lVert#1\right\rVert}

\newcommand{\Cd}{\mathrm{Cauchy}}

\begin{document}
\title{An auditory cortex model for sound processing}
\author{Rand Asswad\inst{1}\and 
Ugo Boscain\inst{2}\orcidID{0000-0001-5450-275X} \and
Giuseppina Turco\inst{3}\orcidID{0000-0002-5963-1857}\and
Dario Prandi\inst{1}\orcidID{0000-0002-8156-5526} \and
Ludovic Sacchelli\inst{4}\orcidID{0000-0003-3838-9448}}
\authorrunning{F. Author et al.}
%
\institute{Université Paris-Saclay, CNRS, CentraleSupèlec, Laboratoire des signaux et systèmes, 91190, Gif-sur-Yvette, France \\ \email{\{rand.asswad, dario.prandi\}@centralesupelec.fr} \and
    CNRS, LJLL, Sorbonne Université, Université de Paris, Inria, Paris, France\\
\email{ugo.boscain@upmc.fr}\\
\and
CNRS, Laboratoire de Linguistique Formelle, UMR 7110, Université de Paris, Paris, France
\email{gturco@linguist.univ-paris-diderot.fr} \\
\and
Université Lyon, Université Claude Bernard Lyon 1, CNRS, LAGEPP UMR 5007, 43 bd du 11 novembre 1918, F-69100, Villeurbanne, France\\
\email{ludovic.sacchelli@univ-lyon1.fr}
}
\maketitle              
\begin{abstract}
     The reconstruction mechanisms built by the human auditory system during sound reconstruction are still a matter of debate. The purpose of this study is to refine the auditory cortex model introduced in \cite{boscainBioinspired2021}, and inspired by the geometrical modelling of vision. The algorithm transforms the degraded sound in an 'image' in the time-frequency domain via a short-time Fourier transform. Such an image is then lifted in the Heisenberg group and it is reconstructed via a Wilson-Cowan differo-integral equation. Numerical experiments on a library of speech recordings are provided, showing the good reconstruction properties of the algorithm. 

\end{abstract}

\section{Introduction}

Human capacity for speech recognition with reduced intelligibility has mostly been studied from a phenomenological and descriptive point of view (see \cite{Mattys} for a review on noise in speech, as well as a wide range of situations in \cite{Assmann2004,Luce}). 
What is lacking is a proper mathematical model informing us on how the human auditory system is able to reconstruct a degraded speech sound.
The aim of this study is to provide a neuro-geometric model for sound reconstruction based on the description of the functional architecture of the auditory cortex.

Knowledge on the functional architecture of the auditory cortex and the principles of auditory perception are limited. For that reason, we turn to recent advances in the mathematical modeling of the functional architecture of the primary visual cortex and the processing of visual inputs \cite{Petitot,Citti2006,boscain2010anthropomorphic} (which recently yield very successful applications to image processing \cite{Franken09,Prandi2017,remizovJMIV}
) to extrapolate a model of the auditory cortex.
This idea is not new: neuroscientists take models of V1 as a starting point for understanding the auditory system (see, e.g., \cite{Nelken2011}). Indeed, biological similarities between the structure of the primary visual cortex (V1) and the primary auditory cortex (A1) are well-known to exist. 
V1 and A1 share a “topographic” organization, a general principle determining how visual and auditory inputs are mapped to those neurons responsible for their processing \cite{Rauschecker}. Furthermore, the existence of receptive fields of neurons in V1 and A1 that allow for a subdivision of neurons in “simple” and “complex” cells supports the idea of a “common canonical processing algorithm within cortical columns” \cite{Tian7892}.  

\section{The contact space approach in V1}
The neuro-geometric model of V1 finds its roots in the experimental results of Hubel and Wiesel \cite{hubel-wiesel}. 
This gave rise to the so-called sub-Riemannian model of V1 in \cite{Petitot1999,Citti2006,boscain2010anthropomorphic,Prandi2017}. 
The main idea behind this model is that an image, seen as a function $f:\R^2\to\R_+$ representing the grey level, is lifted to a distribution on $\R^2\times P^1$, the bundle of directions of the plane.
Here, $P^1$ is the projective line, i.e., $P^1 = \R/\pi\mathbb{Z}$. 
More precisely, the lift is given by $Lf(x,y,\theta) = \delta_{Sf}(x,y,\theta)f(x,y)$ where $\delta_{S_f}$ is the Dirac mass supported on the set $S_f\subset \R^2\times P^1$ of points $(x,y,\theta)$ such that $\theta$ is the direction of the tangent line to $f$ at $(x,y)$.

When $f$ is corrupted (i.e. when $f$ is not defined in some region of the plane), the reconstruction is obtained by applying a deeply anisotropic diffusion mimicking the flow of information along the horizontal and vertical connections of V1, with initial condition $L_f$. This diffusion is known as  the {\em sub-Riemannian  diffusion} in  $\R^2\times P^1$, cf.\ \cite{Agrachev2019}.
One of the main features of this diffusion is that it is invariant by rototranslation of the plane, a feature that will not be possible to translate to the case of sounds, due to the special role of the time variable.

The V1-inspired pipeline is then the following: first a lift of the input signal to an adequate contact space, then a processing of the lifted signal according to sub-Riemannian diffusions, then projection of the output to the signal space.

\section{The model of A1}

\begin{wrapfigure}{r}{.3\textwidth}
    \centering
    \scalebox{.25}{
        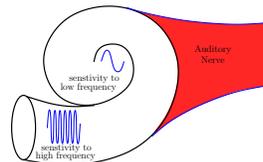
    }
\caption{Perceived pitch of a sound depends on the stimulated location in the cochlea.
}
  \label{fig:cochlea}
\end{wrapfigure}

The sensory input reaching A1 comes directly from the cochlea. The sensors are tonotopically organized (in a frequency-specific fashion), with cells close to the base of the ganglion being more sensitive to low-frequency sounds and cells near the apex more sensitive to high-frequency sounds, see Figure~\ref{fig:cochlea}. This implies that sound is transmitted to the primary auditory cortex A1 in the form of a ‘spectrogram’: when  a sound $s:[0,T]\to\R$ is heard,  A1 is fed with its time-frequency representation
$S:[0,T]\times\R\to \C$. If $s\in L^2(\R^2)$, as given by the short-time Fourier transform of $s$, that is
\begin{equation}
  S(\tau, \omega) := \operatorname{STFT}(s)(\tau,\omega)= \int_{\R} s(t)W(\tau-t)e^{2\pi i t\omega}\,dt.
\end{equation}
Here, $W:\R\to [0,1]$ is a compactly supported (smooth) window, so that $S\in L^2(\R^2)$.
The function $S$ depends on two variables: the first one is time, that here we indicate with the letter $\tau$, and the second one is frequency, denoted by $\omega$. 
Since $S$ is complex-valued, it can be thought as the collection of two black-and-white images: $|S|$ and $\arg S$.
Roughly speaking, $|S(\tau,\omega)|$ represents the strength of the presence of the frequency $\omega$ at time $\tau$. In the following, we call $S$ the sound image. The V1-inspired approach is then to apply image perception models on $|S|$. However, in this case time plays a special role: the whole sound image does not reach the auditory cortex simultaneously, but sequentially, and invariance by image rotations is lost.
As a consequence, different symmetries have to be taken into account and a different model for both the lift and the processing in the contact space are required.


\subsection{The lift procedure}

We propose an extension of the time-frequency representation of a sound, which is at the core of the proposed algorithm. Sensitivity to variational information, such as the tangent to curves for V1, is now translated in a manner that takes into account the role of time. A regular curve $t\mapsto(t, \freq(t))$ in the time frequency domain $(\tau,\omega)$ is lifted to a 3-dimensional \emph{augmented  space} by adding a new variable $\slope=d\freq/d\tau$. The variation $\slope$ in instantaneous frequency is then associated to the chirpiness of the sound.
Alternatively, a curve 
$t\mapsto(\tau(t),\omega(t),\nu(t))$ is a lift of planar curve $t\mapsto (t,\omega(t))$ if $\tau(t)=t$  and if $\nu(t)=d\omega/dt$. Setting $u(t)=d\nu/dt$  we can say that a curve in the contact space $t\mapsto(\tau(t),\omega(t),\nu(t))$ is a lift of a planar curve if there exists a function $u(t)$   such that:
\begin{equation}\label{eq:heis}
  \frac{d}{dt}
 \begin{pmatrix}
    \tau \\ \freq \\ \acc 
 \end{pmatrix}
  =
	\begin{pmatrix}
	    1 \\ \acc \\ 0
\end{pmatrix}
  +u
	\begin{pmatrix}
	0 \\ 0 \\ 1
	\end{pmatrix}=X_0(\tau,\freq,\acc)+u X_1(\tau,\freq,\acc)
\end{equation}
The vector fields $(X_0,X_1)$ generate the Heisenberg group. However, we are not dealing here with the proper sub-Riemannian distribution, since $\{X_0+uX_1\mid u\in\R\}$ is only one-dimensional.

Following \cite{boscain2010anthropomorphic}, when $s$ is a general sound signal, we lift each level line of $|S|$. By the implicit function theorem, this yields the following subset of the augmented space:
\begin{equation}\label{eq:Sigma}
    \Sigma = 
    \left\{ (\tau,\freq,\slope)\in\R^3\mid \slope \partial_\freq |S|(\tau,\freq)+\partial_\tau|S|(\tau,\freq) =0\right\}.
\end{equation}
The external input from the cochlea to the augmented space is then given by
\begin{equation}\label{eq:lift-distr}
    I(\tau,\freq,\slope) 
    = 
    S(\tau,\freq)\delta_{\Sigma}(\tau,\freq,\slope)
    =
    \begin{cases}
        S(\tau,\freq)& \text{if }\nu\partial_\freq|S|(\tau,\freq) =-\partial_\tau|S|(\tau,\freq),\\
        0 &\text{otherwise},
    \end{cases}
\end{equation}
with $\delta_\Sigma$ denoting the Dirac delta distribution concentrated on $\Sigma$.

\subsection{Associated interaction kernel.}
Considering A1 as a slice of the augmented space allows to deduce a natural structure for neuron connections. For a single time-varying frequency $t\mapsto\omega(t)$, its lift is concentrated on the curve $t\mapsto (\omega(t),\nu(t))$, such that
\begin{equation}\label{eq:contr-sys}
    \frac{d}{dt}\begin{pmatrix}
         \omega  \\
         \nu
    \end{pmatrix} 
    =
    Y_0(\omega,\nu) + u(t) Y_1(\omega,\nu),
\end{equation}
where $Y_0(\omega,\nu) = (\nu,0)^\top$, $Y_1(\omega,\nu) = (0,1)^\top$, and $u:[0,T]\to \R$.

As in the case of V1 \cite{Remizov2013}, we model neuronal connections via these dynamics.
In practice, this amounts to assume that the excitation starting at a neuron $X_0=(\freq',\acc')$ evolves as the stochastic process $\{A_t\}_{t\ge 0}$ following the SDE $dA_t = Y_0(A_t) dt + Y_1(A_t) dW_T$, $\{W_t\}_{t\ge 0}$ being a Wiener process, with inital condition $A_0 = (\freq',\acc')$.
The generator of $\{A_t\}_{t\ge 0}$ is the operator $\mathcal L = Y_0 + (Y_1)^2$.

The influence $k_\delta(\freq,\acc\|\freq',\acc')$ of neuron $(\freq',\acc' )$ on neuron $(\freq,\acc)$ at time $\delta>0$ is then modeled by the transition density of the process $\{A_t\}_{t\ge0}$, given by the integral kernel at time $\delta$ of the Fokker-Planck equation 
\begin{equation}
    \label{eq:kolmo}
\partial_t I = \mathcal L^* I , 
\quad\text{where}\quad
\mathcal L^* =- Y_0+(Y_1)^2 = -\acc \partial_\freq + b\partial_\acc^2.
\end{equation}
The vector fields $Y_0$ and $Y_1$ are interpreted as first-order differential operators and the scaling parameter $b>0$ models the relative strength of the two terms. 
We stress that an explicit expression of $k_\delta$ is well-known (see, for instance \cite{Barilari2017}) and is a consequence of the hypoellipticity of $(\partial_t-\mathcal L^*)$.

\subsection{Processing of the lifted signal.}
As we already mentioned, the special role played by time in sound signals does not permit to model the flow of information as a pure hypoelliptic diffusion, as was done for static images in V1. We thus turn to a different kind of model: Wilson-Cowan integro-differential equations \cite{Wilson1972}. This model has been successfully applied to describe the evolution of neural activations, in V1 in particular, where it allowed to predict complex perceptual phenomena such as the emergence of patterns of hallucinatory \cite{Ermentrout,Bressloff01} or illusory \cite{SSVM2019,JNP2019,BertalmioCalatroniEtAl2019} nature.  It has also been used in various computational models of the auditory cortex \cite{Loebel2007,Rankin2015,Zulfiqar2020}.

Wilson-Cowan equations also present many advantages from the point of view of A1 modelling: i) they can be applied independently of the underlying structure, which is only encoded in the kernel of the integral term; ii) they allow for a natural implementation of delay terms in the interactions; iii) they can be easily tuned via few parameters with a clear effect on the results.

On the basis of these positive results, we emulate this approach in the A1 context.
Namely, we consider the lifted sound image $I(\tau,\freq,\slope)$ to yield an A1 activation 
 $a:[0,T]\times \R\times \R\to \C$, solution of 
\begin{multline}\label{eq:WC}\tag{WC}
        \partial_t a(t,\freq,\acc) = {-\alpha a(t,\freq,\acc)}+{\beta I(t,\freq,\acc)}\\
        +{\gamma \int_{\mathbb R^2} k_\delta(\freq,\acc\|\freq',\acc')\sigma(a(t-\delta, \freq',\acc'))\,d\freq'\,d\acc'},
\end{multline}
with initial condition $a(t,\cdot,\cdot)\equiv 0$ for $t\le 0$.
Here, $\alpha,\beta, \gamma>0$ are parameters, $k_\delta$ is the interaction kernel, and $\sigma:\C\to \C$ is a (non-linear) saturation function, or sigmoid. We pick $\sigma(\rho e^{i\theta})=\tilde\sigma(\rho) e^{i\theta}$ where $\tilde \sigma(x) = \min\{1,\max\{0,\kappa x\}\}$, $x\in\R$, for some fixed $\kappa>0$. 

The presence of a delay $\delta$ in the activation appearing in the integral interaction term \eqref{eq:WC}  models the fact that the time-scale of the input signal and of the neuronal activation are comparable.
When $\gamma = 0$, equation \eqref{eq:WC} is a standard low-pass filter $\partial_t a = -\alpha a + I$.
Setting $\gamma\neq 0$ adds a non-linear delayed interaction term on top of this exponential smoothing, encoding the inhibitory and excitatory interconnections between neurons.

\subsection{Algorithm pipeline.} 
Both operations in the lift procedure are invertible: the STFT by inverse STFT, and the lift by integration along the $\acc$ variable (that is, summation of the discretized solution). The final output signal is thus obtained by applying the inverse of the pre-processing (integration then inverse STFT) to the solution $a$ of \eqref{eq:WC}. That is, the resulting signal is given by
\begin{equation}
  \hat s(t) = \operatorname{STFT}^{-1}
\left(
\int_{-\infty}^{+\infty} a(t,\freq,\acc)\diff \acc
\right).
\end{equation}
It is guaranteed that $\hat s$ is real-valued and thus correctly represents a sound signal. From the numerical point of view, this implies that we can focus on solutions of \eqref{eq:WC} in the half-space $\{\freq\geq 0\}$, which can then be extended to the whole space by mirror symmetry.

The resulting algorithm to process a sound signal $s:[0,T]\to \R$ is then:
\begin{enumerate}
  \item[A.] \textbf{Preprocessing:} 
  \begin{enumerate}
    \item Compute the time-frequency representation $S:[0,T]\times \R\to \C$ of $s$, via standard short time Fourier transform (STFT);
    \item Lift this representation to the Heisenberg group, which encodes redundant information about chirpiness, obtaining $I:[0,T]\times\R\times\R\to \C$;
  \end{enumerate}
  \item[B.] \textbf{Processing:} Process the lifted representation $I$ via a Wilson-Cowan equation adapted to the Heisenberg structure, obtaining $a:[0,T]\times\R\times\R\to \C$.
  \item[C.] \textbf{Postprocessing:} Project $a$ to the processed time-frequency representation $\hat S:[0,T]\times \mathbb R\to \mathbb C$ and then apply an inverse STFT to obtain the resulting sound signal $\hat s:[0,T]\to \R$.
\end{enumerate}

\section{Numerical implementation}

The discretization of the time and frequency domains is determined by the sampling rate of the original signal and the window size chosen in the STFT procedure. 
That is, by the Nyquist-Shannon sampling theorem, for a temporal sampling rate $\delta t$ and a window size of $T_w$, we consider the frequencies $\omega$ such that $|\omega|<1/(2\delta t)$, with a finest discretization rate of $1/(2T_w)$.
Observe, in particular, that the frequency domain is bounded.
Nevertheless, the chirpiness $\nu$ defined as $\nu\partial_\omega\abs{S}(\tau,\omega) + \partial_\tau\abs{S}(\tau,\omega)=0$ is unbounded and, since generically there exists points such that $\partial_\omega |S| (\tau_0,\omega_0)=0$, it stretches over the entire real line.

\begin{figure}[t]
    \begin{minipage}{.47\linewidth}
    \centering
    \includegraphics[width=\textwidth]{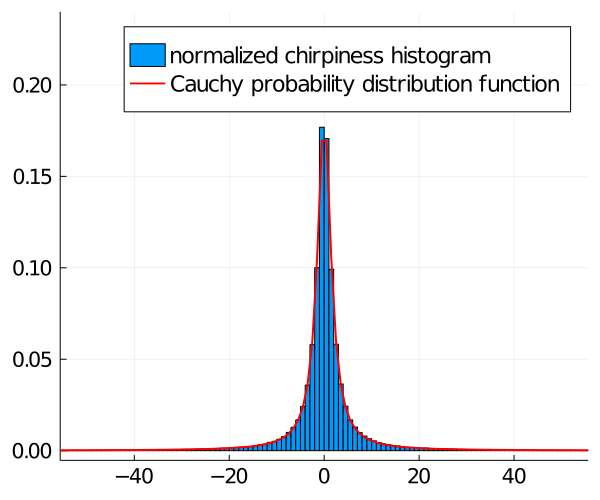}
    \caption{Chirpiness of a speech signal compared to Cauchy distribution}
    \label{fig:cauchy_dist}
    \end{minipage}%
    \hfill
    \begin{minipage}{.47\linewidth}
    \centering
        \includegraphics[width=.27\textwidth]{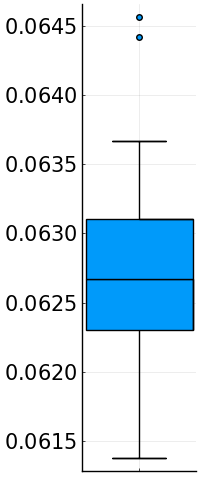}\hspace{0.1\linewidth}%
        \includegraphics[width=.27\textwidth]{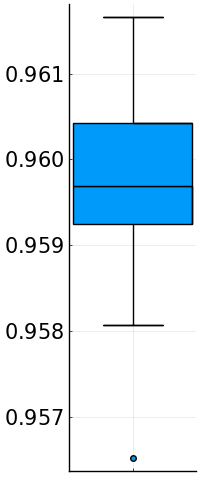}
    \caption{Box plots for estimated Cauchy distributions of speech signals
    chirpiness. \emph{Left:} Kolmogorov-Smirnov statistic values. \emph{Right:} percentage of values falling in $I_{0.95}$}
    \end{minipage}
\end{figure}

To overcome this problem, a natural strategy is to model the chirpiness values as
a random variable, and considering only chirpinesses falling inside the confidence
interval $I_p$ for some reasonable $p$-value (e.g., $p=0.95$). 
A reasonable assumption for the distribution of $X$ is that it follows
a Cauchy distribution $\Cd(x_0, \gamma)$.
Indeed, this corresponds to assuming that $\partial_\omega |S|$ and
$\partial_\tau |S|$ are normal (and independent) distributions \cite{papoulis1991}.
As it is customary, we chose as estimator for the location parameter $x_0$ the median of $X$
and for the scale parameter $\gamma$ half the interquartile distance.

Although statistical tests on a library of real-world speech signals 
\footnotemark[\value{footnote}]
rejected the assumption that $X\sim \Cd(x_0,\gamma)$,
the fit is quite good according to the Kolmogorov-Smirnov statistic
$D_n=\sup_x\abs{F_n(x)-F_X(x)}$. Here, $F_X$ is the cumulative distribution function
of $X$ and $F_n$ is the empirical distribution function
evaluated over the chirpiness values.

\footnotetext{The speech material used in the current study is part of an ongoing psycholinguistic project on spoken word recognition. Speech material comprises 49 Italian words and 118 French words. The two sets of words were produced by two (40-year-old) female speakers (a French monolingual speaker and an Italian monolingual speaker) and recorded using a headset microphone AKG C 410 and a Roland Quad Capture audio interface.
Recordings took place in the soundproof cabin of the Laboratoire de Phonétique et Phonologie (LPP) of Université de Paris Sorbonne-Nouvelle.
Both informants were told to read the set of words as fluently and naturally as possible.}

\section{Denoising experiments}

For simple experiments on synthetic sounds, highligting the characteristics of the proposed algorithm, we refer to \cite{boscainBioinspired2021}.

\begin{figure}[t]
    \begin{minipage}{.47\linewidth}
    \centering
    \includegraphics[width=\textwidth]{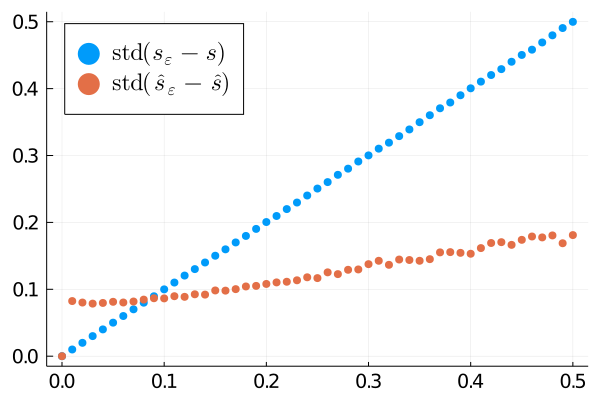}
    \end{minipage}%
    \hfill
    \begin{minipage}{.47\linewidth}
    \centering
    \includegraphics[width=\textwidth]{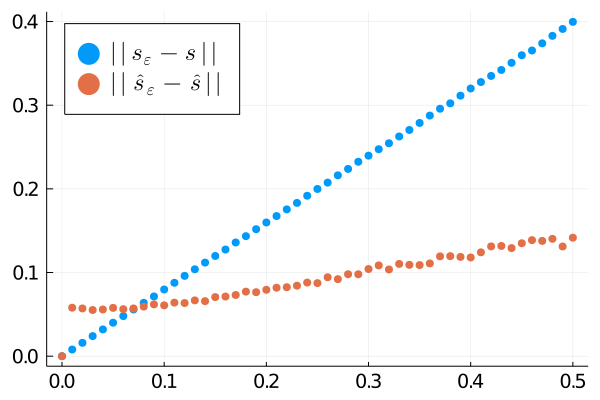}
    \end{minipage}
    \caption{\label{fig:experiments}Distance of noisy sound to original one before (blue) and after (red) the processing, plotted against the standard deviation of the noise ($\varepsilon$). \emph{Left:} standard deviation metric. \emph{Right:} $\norm{\cdot}$ norm.}
\end{figure}

In Figure~\ref{fig:experiments} we present the results of the algorithm applied to a denoising task.
Namely, given a sound signal $s$, we let $s_\varepsilon = s + g_\varepsilon$,
where $g_\varepsilon\sim \mathcal N(0,\varepsilon)$ is a gaussian random variable.
We then apply the proposed sound processing algorithm to obtain $\widehat{s_\varepsilon}$.
As a reconstruction metric we present both the norm $\norm{\cdot}$
where for a real signal $s$, $\norm{s} = \norm{s}_1/\dim(s)$
with $\norm{\cdot}_1$ as the $L_1$ norm
and the standard deviation $\mathrm{std}(\widehat{s_\varepsilon} -\widehat s)$.
We observe that according to both metrics the algorithm indeed improves the signal.

\bibliographystyle{splncs04}
\bibliography{biblio}
\end{document}